\documentstyle[12pt]{article}
\input epsf.sty
\def\thebibliography#1{\section*{{\normalsize \begin{center} \bf
References \end{center} }}\list
  {[\arabic{enumi}]}{\settowidth\labelwidth{[#1]}\leftmargin\labelwidth
    \advance\leftmargin\labelsep
    \usecounter{enumi}}
    \def\newblock{\hskip .11em plus .33em minus -.07em}
    \sloppy
    \sfcode`\.=1000\relax}

\newcommand{\psdiag}[3]{\hspace{1mm}\raisebox{-#1mm}{\epsfysize#2mm
\epsffile{#3.eps}}\hspace{1mm}}

\newtheorem{th}{Theorem}

\newtheorem{lem}[th]{Lemma}

\newcommand{\nc}{\newcommand}

\newtheorem{definition}[th]{Definition}
\newtheorem{prop1}[th]{Proposition}
\newtheorem{lemma}[th]{Lemma}
\newtheorem{remark}[th]{Remark}
\newtheorem{cor}[th]{Corollary}
\newtheorem{example}[th]{Example}

\nc{\bother}{\begin{other}}
\nc{\eother}{\end{other}}
\nc{\bex}{\begin{example}}
\nc{\eex}{\end{example}}
\def\be{\begin{equation}}
\def\ee{\end{equation}}
\nc{\bea}{\begin{eqnarray}}
\nc{\eea}{\end{eqnarray}}
\nc{\bth}{\begin{th}}
\nc{\eth}{\end{th}}
\nc{\bdf}{\begin{definition}}
\nc{\edf}{\end{definition}}
\nc{\bpr}{\begin{prop1}}
\nc{\epr}{\end{prop1}}
\nc{\blm}{\begin{lemma}}
\nc{\elm}{\end{lemma}}
\nc{\br }{\begin{remark}}
\nc{\er}{\end{remark}}
\nc{\bc}{\begin{center}}
\nc{\ec}{\end{center}}

\def\v8{\vskip 8pt}

\def\H{{\cal H}}
\def\s{\sigma}
\def\Z{{\Bbb Z}}
\def\build#1_#2^#3{\mathrel{\mathop{\kern 0pt#1}\limits_{#2}^{#3}}}

\def\l{\lambda}
\def\a{\alpha}
\def\b{\beta}
\def\s{\sigma}

\def\var{\varepsilon}

\newcommand{\CUP}{\mathop{\cup}}

\date{\today}

\input epsf.sty

\begin{document}

\input amssym.def
\input amssym.tex

\begin{center}
{ \LARGE \sl
 Refined  invariants and TQFT's\\[0.3cm] from Homfly skein theory}

\footnotetext{Date:  June 1998}
\v8
\v8
\v8
{\sl Anna Beliakova} \footnotemark[1]

\footnotetext[1]{supported by the Swiss National Science Foundation}

\v8\v8
%

\v8

\parbox{10cm}{{\small {\bf Abstract}:  We work in the reduced $SU(N,K)$ 
modular category as constructed recently by Blanchet.
We define  spin type and cohomological
refinements of the Turaev-Viro  invariants of closed 
oriented 3-manifolds and  give a
formula relating them to  Blanchet's invariants.
Roberts' definition of the Turaev-Viro state sum
 is exploited. Furthermore,
we construct refined Turaev-Viro and Reshetikhin-Turaev TQFT's 
and study the relationship between them. 
}}


 \end{center}
\v8\v8

{\Large \sl Introduction}

\v8\v8
In \cite{t} Turaev
 reduced the construction of quantum 3-manifold 
invariants and TQFT's to the construction of modular categories.
A modular category is a monoidal category with additional structure
(braiding, twist, duality, finite set of simple objects satisfying a 
domination property and a  non-degeneracy axiom). A first example 
of the modular category was obtained from the representation theory of
the quantum group $U_q(sl(2))$. 
Later an elementary  approach, based on the
Kauffman skein relations and
leading to 
the same family of invariants, was developed by Lickorish in [L]. 

Yokota [Y] generalized his approach and constructed
the  $SU(N,K)$ modular category
 using Homfly skein theory. 
The underlying  invariant $\tau^{SU(N)}$  coincides with
the  invariant   of 
 Turaev-Wenzl [TW] extracted from $U_q(sl(N))$ at level $K$.
Recently Blanchet [Bl] defined the
 reduced $SU(N,K)$
modular category.
His
 invariant $\tau$ can be
considered as a generalization of $\tau^{PSU(N)}$
to the case when $N$ and $K$ are not coprime. 
For any closed oriented 3-manifold
$M$ holds (see [Bl])
$$\tau^{SU(N)}(M)\,=\,\tau^{U(1)}(M)\; \tau(M)$$
where $\tau^{U(1)}(M)$ is defined in \cite{o}. Blanchet constructed
cohomological and
spin type  refinements of $\tau(M)$ depending 
on the
so-called spin$^d$ structure  on $M$ with $d={\rm gcd}(N,K)$.
He showed
 that $\tau(M)$ splits  into a sum of refined invariants.


\v8
In this article
we work in the reduced $SU(N,K)$ modular category
as constructed by Blanchet. 
We give  a  definition of the  refinement $Z(M,s,h)$ of the
Turaev-Viro state sum
$Z(M)$  
 depending on the spin$^d$
structure $s$  on a closed 3-manifold $M$ and the
first $\Z/d\, \Z$-cohomology class $h$. We show that
\be \label{2}
Z(M)=
\sum_{s,h} Z(M,s,h)\ee
and prove the relation with Blanchet's invariants
\be
\label{3}
 Z(M,s,h)=\tau(M,s)\;\tau(-M, s+h). \ee
Analogous formulas also hold  for  cohomological refinements.
 The  definition of $Z(M,s,h)$
is given in terms of Roberts' chain-mail link.
It turns out that (\ref{2}) and (\ref{3})
 can be proved  by minor modifications of
Roberts' arguments.  
Nevertheless,
we give a different  proof of (\ref{3})
which generalizes directly to the  TQFT operators.
\v8

 In the last section we  construct
  spin  topological quantum field theories
 (TQFT's)  for type $A$ modular categories. 
 In the $SU(2,K)$ case these theories were studied  
in [BM] and [B]. 
The vector space asssociated to a surface with structure
is defined as (a subspace of) a formal linear 
span of special colorings of some trivalent graph. In contrast
to the unspun (or non-refined)
theory, this  vector space for a non-connected surface 
is not equal to the  tensor product of 
spaces associated with connected components.
We define operators corresponding to  spin 3-cobordisms and prove the
gluing property for them.  
Finally, we construct a weak
 spin TQFT which can be regarded as a `zero graded part' of the spin TQFT.
We  show that the unspun theory is the  sum of
weak spin TQFT's.

Using the same approach,  we  extend Roberts' invariant $Z(M,s,h)$
 to  a refined Turaev-Viro TQFT. Here 
in oder to prove the  gluing axiom  we use an analog
of (\ref{3})  for spin 3-cobordisms.

\v8

{\small\sl
 Acknowledgements:
I wish to thank Christian Blanchet for explaining his work
and for  stimulating e-mail conversations.
My special thank is to  Christof Schmidhuber for
 improving my English.
}


\section{ Definitions}

\noindent
{\bf Homfly skein theory}.
 The  manifolds throughout this paper are 
 compact, smooth and oriented. 
By {\it  links}
we mean {\it isotopy classes of framed links}. 
A framing is a trivialization of the normal bundle. This defines an 
orientation on the link.
In all figures we use the blackboard framing convention.

Let  $M$ be a 3-manifold (possibly with a given set of framed oriented 
points on the boundary).
 We denote by $\H(M)$  the $\Bbb C$-vector space of formal sums of
links in $M$ (and 
framed 
arcs in $M$ that meet $\partial M$ in precisely the given 
set of points)
modulo (isotopy keeping boundary points fixed and)
the Homfly skein relation:
\begin{center}
\mbox{\epsfysize=3cm \epsffile{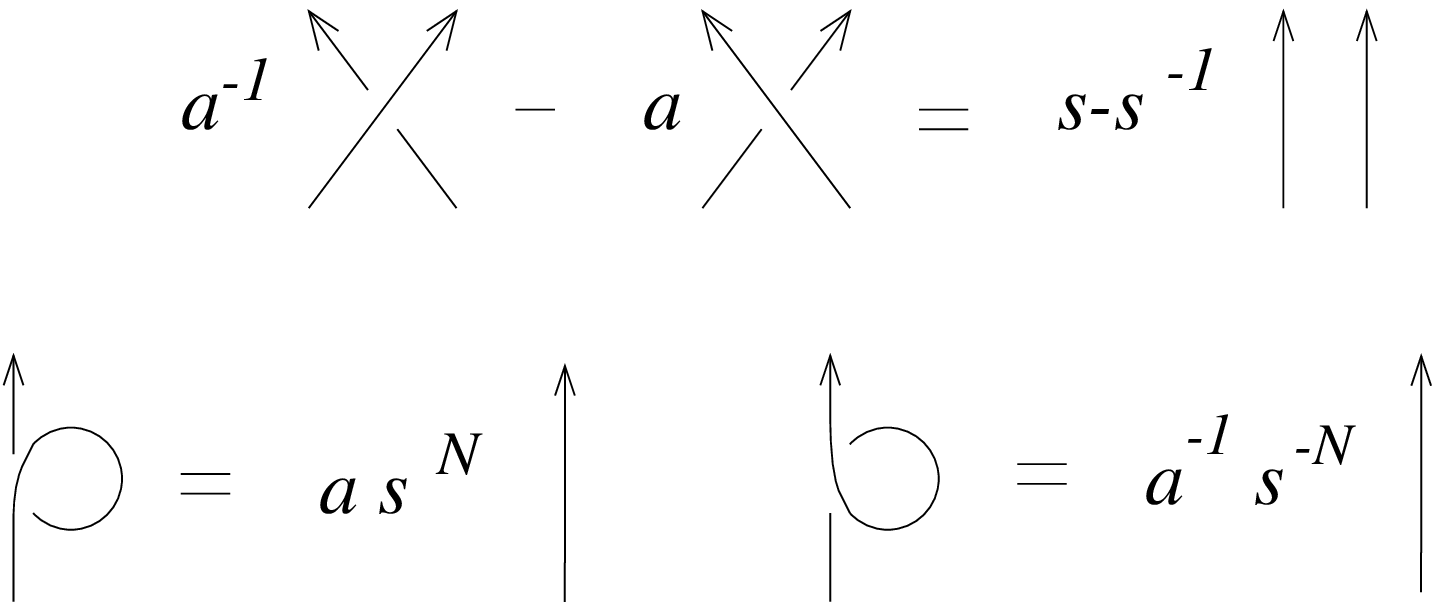}}
\end{center}
$$L\amalg \bigcirc = \frac{s^N-s^{-N}}{s-s^{-1}}L$$
with $a,s\in {\Bbb C}$.
We call $\H(M)$ the skein  of $M$.

For example, $\H(S^3)\cong {\Bbb C}$.
 The isomorphism sends any link $L$
in $S^3$ to its
Homfly polynomial $\langle L\rangle$.

Oriented embeddings induce natural maps between skeins.
Let
\be
\begin{array}{llll}
L_*:& {\H(D^2\times S^1)}^{\otimes m}& \to & \H(S^3)  \\ &
            
x_1\otimes  ...\otimes x_m & \mapsto & \langle L(x_1,..., x_m)\rangle.
\end{array}
\ee
be the map induced by the embedding of $m$ solid tori in $S^3$
with underlying $m$-component link $L$.
We  shall say that the components of $L$ are cabled or colored with $x_1,...,
x_m$.

 \v8
\noindent
{\bf Specification of parameters.}
Let us fix a rank $N\in \Bbb N$ and a level $K \in \Bbb N$, such
that ${\rm gcd}(N,K)=d$ is even,  $N^\prime= N/d$ and
$K^\prime= K/d$ are odd. Let $s$ be a primitive $2(N+K)$ root of unity.
We write $d= \alpha \beta$ with
${\rm gcd}(\alpha, 2K^\prime)={\rm gcd}(\beta, N^\prime)=1$ 
and  choose
 the  framing parameter $a$ such that $(a^N s)^\alpha=1$ and $(a^K s^{-1}
)^\beta =-1$.

\v8\noindent
{\bf Simple objects.}
Denote by $\l=(\l_1,...,\l_p)$ the Young diagram with $\l_i$ boxes in the
$i$-th row. Set $|\l|=\sum^p_{i=1} \l_i$.
In particular, let $1^N$ (resp. $K$) denote the diagram with one
column (resp. one row)  containing $N$ (resp. $K$) cells.

The set of simple objects (colors) in the reduced
$SU(N,K)$ modular category
can be obtained from

$$\{(1^N)^{\otimes i}\otimes \l, \;\;0\leq i<\a,\; \l_1\leq K,\; p<N\}$$
by identifying diagrams which differ by $K^{\otimes \b}$.
Recall
 that for any diagram $\l$ with maximal $N-1$ rows and $K$ columns
$K\otimes \l= K+\l=(K,\l_1,...,\l_p)$.
We denote by $\Gamma_{N,K}$ the resulting set of simple objects.

Under a $\l$-colored line we understand $|\l|$ copies of it
with the  idempotent
of the Hecke algebra corresponding to $\l$ inserted
(see \cite{am} for more details).

 There exists  an involution 
$i:\l\to\l^*$ on the set of colors, such that 
changing  the orientation on the $\l$-colored 
link component is equivalent to 
 replacing  $\l$ by $\l^*$. Note that $|\l|=-|\l^*| \;{\rm mod}\; d$.

\v8\noindent
{\bf Definition of $\omega$.}
Let $y_\l \in \H(D^2\times S^1) $ be the skein element obtained by
cabling with $\l$ a  $0$-framed circle
$\{pt\}\times S^1$, $pt\in D^2-\partial D^2$.
The image of $y_\l$ under the map $\H(D^2\times S^1) \to \H(S^3)$
given  by the standard embedding of the solid torus in $S^3$ is
denoted by $\langle\l\rangle$.

For a cell $c$
in $\l$ with coordinates $(i,j)$ we define its  hook length
 $hl(c)$ and its content $cn(c)$ by formulas
$$ hl(c)= \l_i + \check{\l}_j -i-j+1, \;\;\;\; cn(c)=j-i,$$
where $\check{\l}_j$ is the length of the $j$-th column of $\l$.
Then (see  \cite{a})

\be
\langle\l\rangle=\langle\l^*\rangle=\prod_{cells} \frac{[N+cn(c)]}{[hl(c)]}\;\;\;{\rm where}
\;\;\;[n]=\frac{s^n-s^{-n}}{s-s^{-1}}.\ee

With this notation the element
$$\tilde{\omega}= \sum_{\l\in \Gamma_{N,K}}\langle\l\rangle\; y_\l\;\; \in\;\;
 \H(D^2\times S^1)$$
has the  nice property that the  Homfly polynomial of a
                     link with an   $\tilde{\omega}$-colored 
component  is invariant
under handleslides along this component. In addition, it is also independent 
of the orientation on this component.

We choose the normalization $\omega= \eta \tilde {\omega}$
with
$$\eta^{-2}=\langle \tilde{\omega}\rangle=
(-1)^{\frac{N(N-1)}{2}} \frac{d(N+K)^{N-1}}{\prod^{N-1}_{j=1}
(s^j-s^{-j})^{2(N-j)}}.$$
Then
we have $\langle U_1(\omega)\rangle
\langle U_{-1}(\omega)\rangle=1$ where
$\langle U_\epsilon(\omega)\rangle$ denotes 
   the  
Homfly polynomial on 
the $\epsilon$-framed unknot colored with $\omega$.

\v8
\noindent{\bf Grading.}
The  algebra $\H(D^2\times S^1)$ has a natural  ${\Bbb Z}_{d}=
{\Bbb Z}/d\, {\Bbb Z}$ grading (recall $d={\rm gcd}(N,K)$)
by taking the number of strands modulo $d$.
According to this grading  we decompose
$$\omega=\omega_0 + \omega_1 +... +\omega_{d-1}. $$
Modifying slightly the calculations in Lemma 2.4 \cite{bl}
we get
\be \label{om} \langle
U_0(\omega)\rangle= d\langle
U_0(\omega_i)\rangle= \eta^{-1}.\ee  
By Lemma 4.3 in [Bl] we have $\langle U_1(\omega)\rangle
=\langle U_1(\omega_{d/2})\rangle$.

\v8\noindent
The graded handleslide property can be written as follows
(see \cite{bl}, Lemma 4.1)
\begin{center}
\mbox{\epsfysize=2.5cm\epsffile{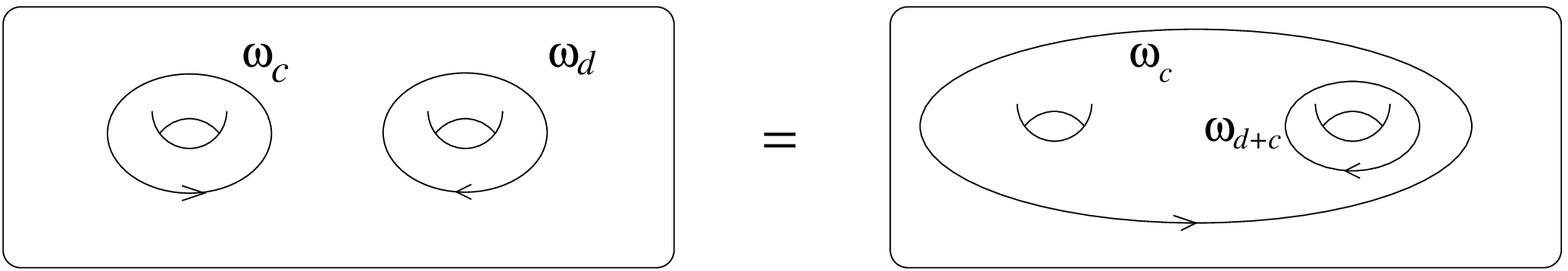}}
\end{center}


\noindent
By Proposition 1.5 in \cite{bl} we have
\begin{center}
\mbox{\epsfysize=1.7cm\epsffile{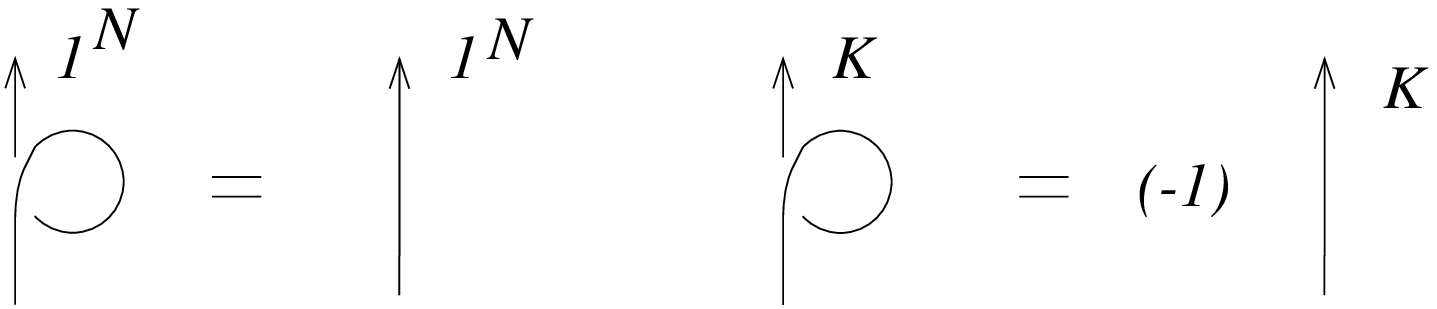}}
\end{center}

\begin{center}
\mbox{\epsfysize=1.7cm\epsffile{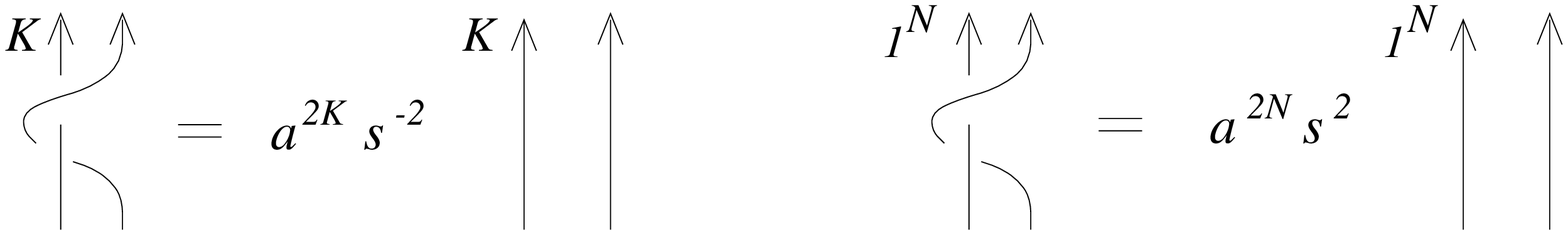}}
\vspace*{0.2cm}
\\ Figure 1: {\it Framing and twisting coefficients on $1^N$ and $K$.}
\end{center}


\v8\noindent 
{\bf Killing property.}
 If $\l\neq 0$,  the following skein element 
\be\label{zero}
\mbox{\epsfysize=1.7cm\epsffile{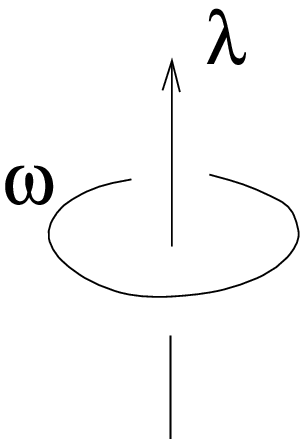}}
\ee
is zero. This is 
an  analog of Lemma 2.5 in [Bl] in  the reduced category.
\v8

\noindent
{\bf Graded killing property.} The skein element (\ref{zero})
with $\omega$ replaced by $\omega_i$ is zero if $\l\neq
(1^N)^{\otimes k} \otimes K^{\otimes l}$
where $0\leq k<\a$ and $0\leq l<\b$
(see  Lemma 4.4 in \cite{bl}).
\v8

\noindent
{\bf Fusion rules.} 
We denote by $\H(D^3, a_1...a_n, b_1...b_m)$ the skein of a 3-ball with 
$n$ outgoing and $m$ incoming
points on the boundary colored with $a_1,...,a_n$ and $b_1,...,b_m$
respectively. Then 
 a natural pairing $\H(D^3, \l\mu, \nu)\times
\H(D^3, \nu, \l\mu)\to \H(S^3)$ can be defined 
 by gluing 3-balls together (identifying
    points of the same color).

With this notation,
the domination property can be written as follows:
\be\label{f}
\sum_{\nu\in \Gamma_{N,K}}\sum_{\alpha}\; \langle \nu\rangle
\psdiag{8}{19}{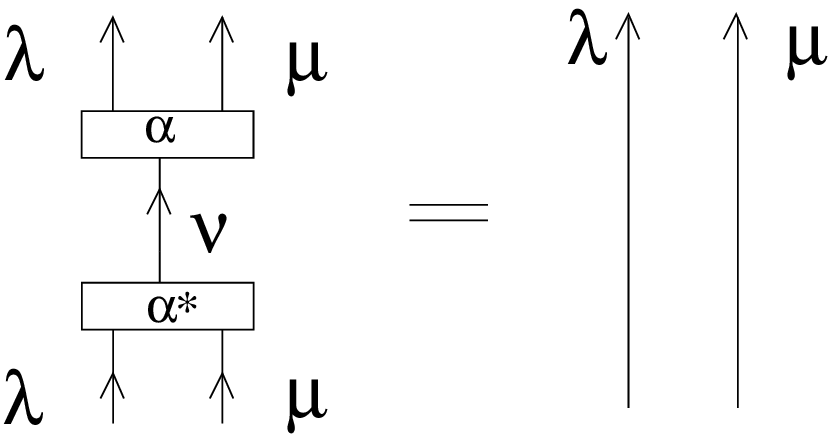},
\ee
where
  $\alpha$ and $\alpha^*$
run over dual  bases 
with respect to the pairing described  above. 
%
%
In what follows
we  shall represent the elements of
  $\H(D^3, \l\mu, \nu)$ by colored
3-vertices for brevity. 
Let $N^\nu_{\l\mu}$ be the dimension
   of  $\H(D^3, \l\mu, \nu)$.
We say that 
a coloring
 $(\l,\mu, \nu)$ of a 3-vertex is admissible if 
 $N^\nu_{\l\mu}\neq 0$. We shall  call  $N^\nu_{\l\mu}$
the  multiplicity of the colored 3-vertex.

We choose a normalization of 3-vertices, so that the following equation hold
(see [BD]):
\be\label{ff}
\psdiag{10}{24}{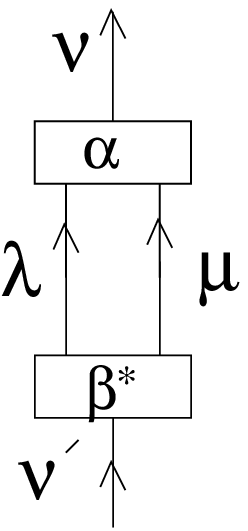}=\delta_{\nu \nu^\prime}\delta_{\alpha \beta}
{\langle \nu\rangle}^{-1}\psdiag{10}{24}{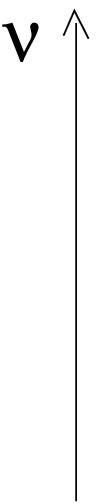}.
\ee
After closing the ends in (\ref{ff}) and applying (\ref{f})
we get $\langle \l\rangle \langle \mu\rangle =\sum_{\nu} N^{\nu}_{\l\mu}
\langle \nu\rangle$.
As a consequence, we have the following rule for deleting
and/or introducing of a $0$-colored line:
\be
\psdiag{5}{15}{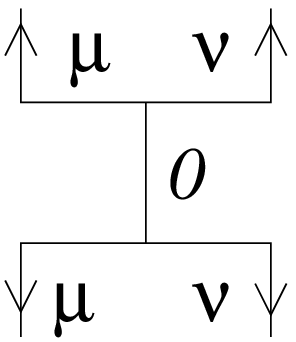}
=\frac{\delta_{\mu^* \nu}}
{\langle\nu\rangle }\psdiag{5}{15}{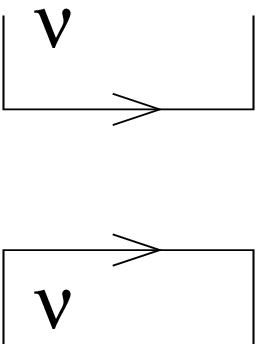}
\ee


\subsection{ Spin$^d$ Structures}
\noindent

{\it
All homology and cohomology groups throughout this article will have 
${\Bbb Z}_d$ coefficients where $d$ is an even integer.}

An oriented  manifold has spin structure if the rotational 
group $SO$ as the
structural group of its stable
tangent bundle can be replaced 
by its 2-fold covering group
 $Spin$ (see [LM, p.80]). 
The notion of a spin$^d$ structure is a natural generalization of this 
construction which
corresponds to the lifting of the structural group
$SO$  to its 
$d$-fold covering group $Spin^d=(Spin\times \Z_d) / \Z_2$  
(where 
${\Bbb Z}_2$ acts by  $(-1,d/2)$ on $Spin\times \Z_d$).
Such structures  always exist on oriented $n$-manifolds 
with $n\leq 3$ due to the 
vanishing of the second Stiefel-Whitney class.

\bdf Let $N$ be an $n$-manifold (possibly
with boundary), where $n=2,3$. Let $FN$ be the space of 
oriented orthonormal $3$-frames on $N$
(= the principle stable tangent bundle). A spin$^d$ structure 
on $N$ is a cohomology class $s\in H^1(FN)$ whose restriction
to each fibre is non-trivial.
\edf

 We denote by Spin$^d(N)$
the set of spin$^d$ structures on $N$.
Using K\"unneth formula one can show that the following  sequence
$$0\to H^1(N)\to H^1(FN)\to H^1(SO(3))\to 0$$
is exact. Therefore,
Spin$^d(N)$
is affinely isomorphic to $H^1(N)$ and  consists of 
$ s\in H^1(FN)$, which are equal to $d/2$ on  homologically trivial 
$0$-framed curves in $N$.


\v8
If a closed 3-manifold
$M=S^3(L)$ is obtained by surgery on $S^3$ along an $m$-component 
link  $L$,
Spin$^d(M)$ is in bijection with the 
 solutions $c=(c_1, ..., c_m)\in ({\Bbb Z_d})^m$ of the following system of
equations
\be\label{char}
\sum^m_{j=1}\; L_{ij}\, c_j = d/2\;  L_{ii}\;\;\;{\rm mod}\;d, \;\;\; \;
1\leq i\leq m\ee
where $\{L_{ij}\}_{1\leq i,j\leq m}$ denotes the linking matrix of $L$
with framing on the diagonal.

\v8
Let $M$ be  a  3-manifold
with parametrized  boundary, i.e. its
 boundary components  are supplied with   diffeomorphisms  
 to the standard  surface. Then we glue (along the parametrization)
to each $\Sigma \in \partial M$  the standard  handlebody.
 The result is a closed 3-manifold $\tilde{M}$.
Deformation retracts of the handlebodies glued to $M$
can be viewed  as 
 a 3-valent  graph  $G$ in $\tilde{M}$ (see Figure 3).
Let $A=\{a_1, ...,a_p\}$ be the set of circles of $G$, where
$p=1-\chi(\partial M)/2$, 
$\chi(\Sigma)$ being  the Euler characteristic of $\Sigma$.
Let $\tilde{M}=S^3(L)$ be obtained by surgery on an $m$-component link $L$. 
Denote by $\{\tilde{L}_{ij}\}$ the linking matrix of $L\cup A$. Then
Spin$^d(M)$ is in bijection
 with the solutions $\tilde{c}=(c_1,...,c_m,z_1,...,z_p)$
of the following equations
\be\label{re}
\sum_j \tilde{L}_{ij} \tilde{c}_j=d/2\; L_{ii} \;\;\;\;{\rm mod}\;d, \;\;\;\; 1\leq i\leq m.
 \ee
The proof in the case of spin structures can be found in [B].
The generalization is straightforward. 
\section{Spin state sum invariants}

{\bf Definition.}
Let $M$ be a closed, connected $3$-manifold.
Choose a handle decomposition of $M$ with $d_0,\, a,\, b,\, d_3$
handles of indices 0,1,2 and 3 respectively.
 Let $H$ be a  handlebody
given by the  union of 0- and 1-handles. Denote by 
$m=\{m_1,...,m_a\}$ and 
 $\var=\{\varepsilon_1, ...,\var_b\} $  
the meridian curves of 1-handles and
the attaching curves of 2-handles on $\partial H$ respectively.
We choose an orientation on all these curves and extract
the  normal vector  
from  the orientation of $\partial H$.
Let  $j(m)$ and $j(\var)$ be  the    images 
 of $m$ and $\var$
under   an orientation preserving embedding $j: H\hookrightarrow S^3$.
Then
$R= j(m)_+\cup j(\varepsilon)_-$
is the  Roberts chain-mail link. Here
  $+$ (resp. $-$)   means
 the  push-off  in the direction 
of the  outgoing  (resp. incoming) normal   to $\partial H$. 

 Let  $s\in {\rm Spin}^d(M)$. Let
$s_0$ be the unique spin$^d$ structure on $S^3$. 
Then $x=s|_H -s_0|_{j(H)}$ assigns $\Z_d$-
numbers  $\{x_1,...,x_{a}\}$
to 1-handles.
Here we assume that the cores of 1-handles are
$0$-framed and  oriented in such a way
that they have the linking number one with the corresponding meridians. 
 Choose  a 2-cycle
$y=\sum_i  y_i \varepsilon_i$
representing a second homology class of $M$. 
Let $h=D(y)\in H^1(M)$
be its Poincare dual class.
We define
$$Z(M,s,h)=(d \eta)^{d_3+d_0-2} \langle
R(\omega_{x_1},...,\omega_{x_{a}}, \omega_{y_1},
...,\omega_{y_{b}})\rangle.$$
Let $-M$ be $M$ with the  reversed orientation. 
By   definition we have  that  $Z(M,s,h)=Z(-M,s,h)$.

\begin{th} $Z(M,s,h)$ is an invariant of $(M,s,h)$.
\end{th}

\noindent
{\bf Proof:} 
 We need to show that $Z(M,s,h)$ does not  depend on
the orientation of $R$,
embedding $j$, the handle decomposition  and 
the representatives for  $x$ and $y$.

Let $\tilde{R}$
 be $R$ with the  orientation on the first component reversed.
After changing the orientation,  
 the numbers $\{-x_1, x_2,..., x_a\}$ will be assigned
to  1-handles.
Applying the
 involution to the set of colors we have
$$\langle R(\omega_{x_1},...)\rangle=\langle \tilde{R}(\omega_{-x_1}, ...)\rangle.$$
 Other cases can be treated analogously.

Two embeddings of $H$ in $S^3$
may be related by unknotting of $1$-handles and 
reframing (twisting of 1-handles across their meridian discs). 

An unknotting move 
can be realized by sliding all 
$\varepsilon$-curves in a  $1$-handle
over a meridian of the other. This does not change the grading
on the meridian, because the boundary of $y$  
is zero and therefore
the number of $\varepsilon$-strands
in  each 1-handle is $0$ modulo $d$.  

Independence of the
reframing move can be shown as follows:
Add to $R$ an $\omega_{d/2}$-colored $\pm 1$-framed 
unknot (unlinked with $R$), slide all $\varepsilon$-curves in  the $i$th
1-handle over it,
twisting them. By the same argument as before
the grading of the unknot remains unchanged.
Finally, slide the unknot over the meridian of this
 1-handle and remove it.
This changes the grading of the meridian by $d/2$, but the coefficient
$x_i$ is also changed by $d/2$ after  reframing.

Two handle decompositions of $M$ can be related by births or deaths of 0-1-,
1-2- and 2-3-handle pairs and handleslides of 1-1- or 2-2-pairs.
The handleslides do not affect the invariant.
Births of 0-1- or 2-3-handle pairs 
add to $R$ a 0-framed unknot which can be slid
 over other 'parallel'
 components and  deleted just like in \cite{r}.
The 1-2-handle pair adds to $R$ a (0,0)-framed Hopf link colored by 
$(\omega_0, \omega_0)$ or a $(\pm 1, 0)$-framed $(\omega_0, 
\omega_{d/2})$-colored one. In both cases the corresponding
 skein elements are equal to one 
by the lemma below.

Representatives of ($x$ or) $y$ differ by changing all labels in the 
(co-)boundary
of  some (0- or) 3-handles. This can be realized by adding a $0$-framed
$\omega_i$-colored unknot, sliding it over all ($m$-curves or)
$\var$-curves in the (co-)boundary  of these handles
and removing it.
$\hfill\Box$
\vspace*{0.2cm}

\noindent{\bf Remark.} By desregarding grading in the proof 
we can see that
$$Z(M)=\eta^{d_0+d_3-2} \langle R(\omega,...,\omega)\rangle$$
is an invariant of $M$.

\begin{lem}\label{dfg}
Let $H_{\epsilon,0}$ be the $(\epsilon,0)$-framed Hopf link with 
$\epsilon=0,\pm1$. 
Then for $i,j \in \Z_d$ we have
\be 
\langle H_{\epsilon,0}(\omega_i, \omega_j)\rangle
=\left\{\begin{array}{ll}1,&{\rm if} 
\;\;\;\epsilon=0, i=0, \;j=0 ;\\
& {\rm or}\;\;\; \epsilon =\pm 1, i=0,\; j=d/2; \\
0,&{\rm otherwise\, .}\end{array}\right. \ee 
\end{lem}

\noindent{\bf Proof:} Graded killing property
 implies that the 
 $\epsilon$-framed component of the Hopf link should be $0$-graded.
Using the identities on Figure 1 we can write
$$\langle H_{\epsilon, 0}(\omega_0,\omega_j)\rangle=\frac{1}{d}
\;\sum^{\a -1}_{k=1} (a^N s)^{2kj}\sum^{\b-1}_{l=1} (-1)^{\epsilon l}
 (a^K s^{-1})^{2lj}$$
which is  non-zero only  in the  two cases mentioned above.
$\hfill\Box$

\begin{th}
For a closed connected 3-manifold $M$, the Turaev-Viro invariant $Z(M)$ 
decomposes as a sum of the refined invariants:
$$Z(M)=\sum_{s,h} Z(M,s,h).$$

\end{th}

\noindent{\bf Proof:} 
The identification of $Z(M)$ with the Turaev-Viro invariant
in the reduced $SU(N,K)$
modular category can be made analogously to
 Theorem 3.6 
in \cite{r1} (see also \cite{it}).  The main difference is that $6j$-symbols
are not numbers, but the elements of the tensor product of four
vector spaces. In the definition of 
the Turaev-Viro state sum a contraction 
over these spaces is  added (see [T] or [BD] for more 
details). 

We will show
the decomposition formula 
in the special case when the
handle decomposition of $M$ is a Heegaard splitting and 
$H$ is embedded standardly in $S^3$. 
For any 
 grading 
of $\var$-curves 
which does not correspond to homology  classes,
 $R$ contains
a meridian curve  linked with $\var$-strands whose total 
grading is not $0$ modulo $d$. This is zero by the  killing property.
If the grading of $m$-curves does not correspond to spin$^d$ 
structures, there exists 
a homologically trivial
 1-cycle in $M$, 
such that the sum over gradings of 1-handles 
representing it is not $0$ mod $d$. After handleslides 
(if necessary) we  represent 
 this cycle by an $\var$-curve.
 Now the invariant vanish by Lemma 4.2 in [Bl].
$\hfill\Box$

\vspace*{0.2cm}

\section{Relation with  Blanchet's  invariants}

In \cite{bl} the refined Reshetikhin-Turaev invariants 
for the reduced $SU(N,K)$ modular category were defined in the following way:
Let $M=S^3(L)$ be given by surgery on $L$. Let $c$ be the solution of
the modulo $d$ characteristic equations (\ref{char})
corresponding to $s\in {\rm Spin}^d(M)$. Then
\be\tau(M,s)=\Delta^{-\sigma(L)}\langle L(\omega_{c_1}, ...,
\omega_{c_m})\rangle
\ee
is  Blanchet's invariant of $(M,s)$. Here $\Delta=
\langle U_1(\omega_{d/2})\rangle$ and
$\sigma(L)$ is the signature of the linking matrix.
This invariant is multiplicative with respect to  connected sums
and normalized at 1 for $S^3$.
We denote  by $\bar{L}$
the mirror of $L$. Then
$$\tau(-M,s)=\Delta^{\sigma(L)}\langle \bar{ L}(\omega_{c_1}, ...,
\omega_{c_m})\rangle . 
$$ 

\begin{th} For a closed connected
3-manifold $M$,
\be \label{23}
 Z(M,s,h)=\tau (M, s+h)\;\tau(-M,s)=\tau(M,s)\;\tau(-M, s+h).\ee
\end{th}

\noindent{\bf Proof:} Once again, we take  Heegaard 
splitting as  handle decomposition
and we choose
the standard embedding of $H$ in $S^3$. Denote by $\Sigma$
the boundary of $H$ with the standard homology basis $\{m_i, 
l_i\}_{1\leq i \leq g}$.
Write $M=H\CUP\limits_\phi -H$ where $\phi:\Sigma\to -\Sigma$ is a gluing 
diffeomorphism.  Note that $\var_i=\phi^{-1}(m_i)$.
 Then  $R= m_+ \cup \var_-$ with grading
$\{x_1, ...,x_g,y_1,...,y_g\}$.
For any link $L$ in a 3-manifold $N$ (possibly with boundary)
we denote by $N(L)$ the result of surgery  on $N$ along $L$.

\v8
Our first aim is to see that
$S^3 (R)=M \# -M$. 
We proceed as follows. Let us cut $S^3$ with $R$ inside 
along $\Sigma$. We get
$$S^3(R)= (S^3-H)(m_+)\;\CUP\; H(\var_-).$$
Once again, cut  out from $H$ a cylinder containing $\var_-$. Then
\be\label{dec} 
S^3(R)= (S^3-H)(m_+)\;\CUP \;(\Sigma\times I)(\var_-)\;\CUP\; H.\ee  
Observe that surgery along $m_+$ on the handlebody
 $S^3-H$   
interchanges the contractible and non-contractible cycles in the 
homology basis of its boundary, i.e
\be\label{13}
\begin{array}{ll}
S^3(R)&= -H\;\CUP\limits_{id}\; (\Sigma\times I)(\var_-)\;
\CUP\limits_{id}\; H =\\  &=
-H\;\CUP\limits_{\phi}\; (\Sigma\times I)(m_-)\;\CUP\limits_{\phi^{-1}}\; H .
\end{array}\ee
Here we have
 used that $\phi(\var_i)=m_i$. Taking into account 
that $(\Sigma\times I)(m_-)$ can be mapped to 
$H \# -H$ by a  diffeomorphism which is the identity on the boundary,
 we get
$$S^3(R)= (-H\CUP\limits_{\phi} H)
 \#(-H \CUP\limits_{\phi^{-1}}  H) =- M \# M.$$  
It remains to find out to which spin$^d$ structure
on $-M\#M$ corresponds 
the grading of $R$. 
According to
the  definition, the spin$^d$ structure on $S^3(R)$ does not 
extend over meridians of not  $0$-graded components  of $R$.
In (\ref{13})
 the structure does not extend over 1-handles of $H$ and $-H$
with $x_i\neq 0$. 
This spin$^d$ structure is equal to
 $s_0+ \sum x_i l_i$ and coincides with 
$s$.
 We have the additional obstruction on $\Sigma\times I$ given by 
 meridians  of curves
${m_i}=\phi(\var_i)$  with $y_i\neq 0$ pushed slightly into interior.
After surgery,  they  become  homologous to $l_i$
on $-H$ and  add 
 the Poincare dual class of $y$ to the spin$^d$
structure on $M$. For
the second equality in (\ref{23}) we use the independence of 
$Z(M,s,h)$   of the orientation of $M$.
$\hfill\Box$

\v8\v8\begin{center}
{\bf Cohomological refinements}
\end{center}
 We need to change the 
specification of parameters in the Homfly polynomial.
The  spin case, considered above,  is here excluded.

For a given rank $N$ and level $K$ choose $s$ be a primitive 
 root of unity of order $2(N+K)$ if $N+K$ is even and 
of order $N+K$ if $N+K$ is odd. As before, $d={\rm gcd}(N,K)$. 
If $N+K$ is even, we suppose that $N^\prime=N/d$ is odd. Then 
$d=\alpha\beta$ with ${\rm gcd}(\alpha, 2K^\prime)=
{\rm gcd}(\beta, N^\prime)=1$ and
we can find  $a$ satisfying
$$(a^Ks^{-1})^\b=(-1)^{N+K+1}, \;\;\;\;\; (a^N s)^\a=1.$$
The main difference to the previous case is that
the twist on the $K$-colored line is trivial and therefore
$\langle U_1(\omega_0)\rangle=\langle U_1(\omega)\rangle=\Delta$.

Let $M=S^3(L)$ and $h\in H^1(M)$. Denote by $c=(c_1, ...,c_m)$
the element in the kernel of the linking matrix (modulo $d$)
corresponding to $h$.
Then 
$$\tau(M,h)=\Delta^{-\s(L)}\langle L(\omega_{c_1},...,\omega_{c_m})\rangle$$
is Blanchet's invariant. Analogously to the spin case, we can define
$Z(M,x,h)$ for any $x\in H^1(M)$ and show its invariance.
The principal  modifications are  that
 the reframing is performed with an $\omega_0$-colored unknot  and that
a birth  of a 1-2-handle pair introduces an
$(\omega_0, \omega_0)$-colored Hopf link
with at least one 0-framed component.  

Analogously we get
$$ Z(M)=\sum_{x,h} Z(M,x,h)\;\;\;\;{\rm
and}\;\;\;\;
 Z(M,x,h)=\tau(M,x)\;\tau(-M, x+h). $$

\section{Spin topological quantum field theories}

A TQFT is a functor from the category of 3-cobordisms to the category 
of finite-dimensional vector spaces. It associates 
 to any closed surface $\Sigma$ 
a vector space $V_\Sigma$  and  to any 3-cobordism $M$
with $\partial M=-\partial_-M\cup\partial_+ M$ an operator
$Z(M):V_{\partial_-M}\to V_{\partial_+ M}$.
Crucial is the 
functorial behavior with respect to the composition of cobordisms 
(gluing property). 
 Two well-known examples of such a construction
are  the Reshetikhin-Turaev (RT) and Turaev-Viro (TV)
 TQFT's (see [T]).



A spin TQFT is a TQFT based on the
 category   $\cal S$  of spin 3-cobordisms.
To define $\cal S$ we need a homotopy-theoretical 
 definition of the notion of a
spin$^d$ structure. 

\bdf
Let $\pi$ be the fibration $BSpin^d\to BSO$. Let $N$ be 
an n-dimensional 
manifold, possibly with boundary.
A $w_2$-structure 
on $N$ is a map $f:N\to BSpin^d$, such that $\pi\circ f$ classifies
the stable tangent bundle of $N$.
A spin$^d$ structure on $N$ is a homotopy class of  $w_2$-structures.

\noindent
Let us fix a  $w_2$-structure on a subset $A\subset N$. 
A relative spin$^d$ structure on $N$ is a homotopy class 
(relative to $A$) of  $w_2$-structures on $N$ extending the 
one given  on $A$.

\edf

 \noindent
{\bf  The category  of spin cobordisms}.
Let $\Sigma$ be a closed
surface. Let us mark a point in each connected 
component of $\Sigma$ and  denote by $P$  the resulting set
of points. We  choose  a $w_2$-structure on $P$.
  Let $\sigma$ be the relative spin$^d$ structure
on $\Sigma$ extending the one given on $P$.
The set of such  structures 
 is affinely isomorphic to
$H^1(\Sigma,P)\cong H^1(\Sigma)$
 by the obstruction theory (see [Sp, p.434]).

 The triple $(\Sigma, P, \sigma)$ is 
an object of $\cal S$.
A morphism
from $(\Sigma,P, \sigma)$
to $(\Sigma^\prime,P^\prime, \sigma^\prime)$ is a 3-cobordism
$M$ with $\partial M=-\Sigma\amalg\Sigma^\prime$ supplied with 
a relative spin$^d$ structure extending the one on $P\cup P^\prime$,
such that its  restriction to the boundary is equal to
 $\sigma\amalg\sigma^\prime$.
The set of such  structures 
on $M$ is affinely isomorphic to
$H^1(M,\partial M)$ (use the exact  sequence
$0\to H^1(M,\partial M)\to H^1(M,P\cup P^\prime)\to H^1(\partial M, P\cup
P^\prime)\to ...$). Here we identify $ H^1(M,\partial M)$ with its image in 
$H^1(M,P\cup P^\prime)$.
 
Let us assume that the boundary of $M$ is parametrized.
Then we can extend the
parametrization diffeomorphism to the map $
M\to \partial M$,
which composed with  $\sigma\amalg \sigma^\prime$
 defines the relative spin$^d$ structure $\tilde{\sigma}$ on $M$.
Any other relative spin$^d$ structure on $M$ 
(with the given restriction to the boundary) is of the form
 $\tilde{\sigma} + 
 \dot{s}$ for some $\dot{s}\in H^1(M,\partial M)$.


\v8
\noindent
{\bf  Spin RT TQFT.}
Let $(\Sigma,P, \sigma) \in {\rm Ob}({\cal S})$
consist  of $n$ connected components. 
Let $\phi: \Sigma\to \Sigma^{st}$ be the parametrization diffeomorphism
respecting the order of components and
$\Sigma^{st}=\Sigma_{g_1}\cup...\cup\Sigma_{g_n}$.
Let us construct a 
framed graph
$\hat{G}^\Sigma$ by taking the  graph $\hat{G}^{g_1}\cup...
\cup \hat{G}^{g_n} $ 
(see Figure 2)   and by
connecting its 1-vertices 
with a fixed  trivalent graph $F_n$. 
\begin{center}
\mbox{\epsfysize=1.7cm\epsffile{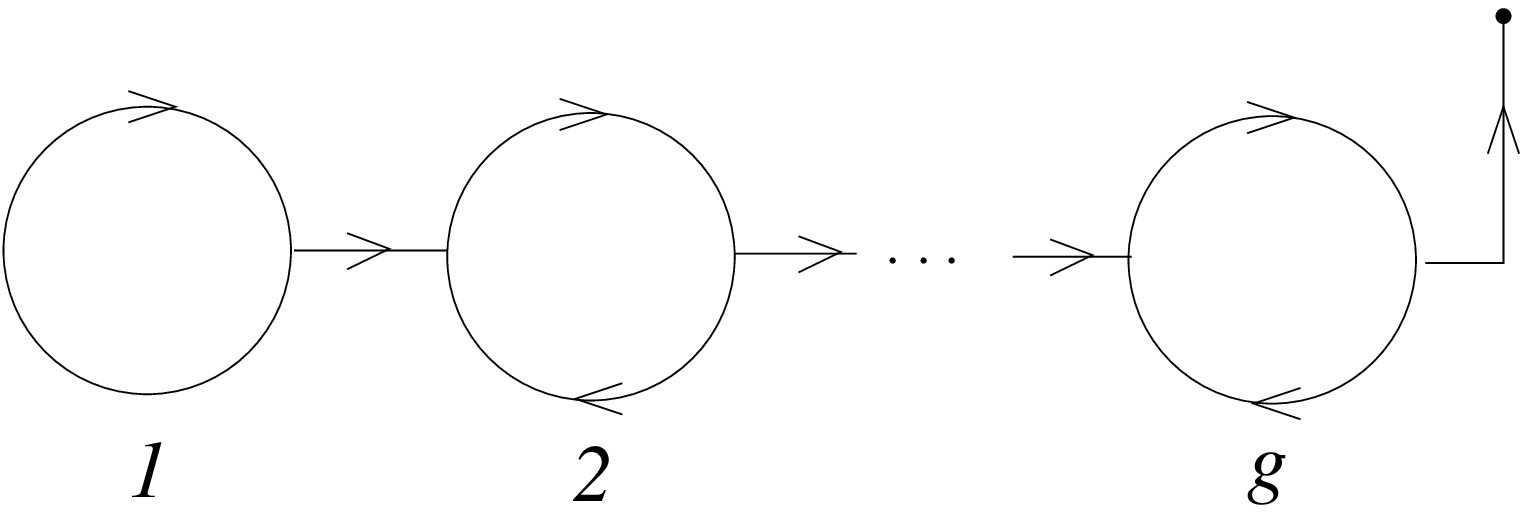}}
\\ Figure 2 {\it The graph $\hat{G}^g$} 
\end{center}

As before, we denote by $\{m_i\}$ the $0$-framed
meridians of the standardly
embedded surface $\Sigma^{\rm st}$. Let $z_i$ be the result of 
the evaluation of the cohomology class corresponding to $\sigma$
on the homology class of  $\phi(m_i)$.

Under a {\it special coloring} $e$ of $\hat{G}^\Sigma$ we understand an
admissible coloring of $\hat{G}^\Sigma$, such that  
 the grading of  colors on the  $i$th
circle is equal to $z_i$\footnotemark[2],
 and  the color of the $i$th line of $F_n$  is 
$(1^N)^{\otimes k_i}\otimes K^{l_i}$ with
$0\leq k_i < \alpha$ and $0\leq l_i<\beta$. 
We denote by $\hat{G}^\Sigma_e$ the  $e$-colored graph.
\footnotetext[2]{Note that the grading is well-defined on the circles of 
$\hat{G}$,
because all  lines connecting two circles are  $0$-graded.}
We  set
$\langle e\rangle=\prod_{e_i\in e}
\langle e_i\rangle$ if ${\rm card}( e)>1$ and $\langle e\rangle=1$ otherwise. 

Let
us  choose the numbering of the lines
of $F_n$, so that the line containing
$k$th 1-vertex
becomes  the number $k$. Then for an ordered set
 $u=\{0,u_2,...,u_{n},0...,0\}$ of $2n-3$ elements
we define
$u_e=(a^N s)^{2\sum k_i u_i} (a^K s^{-1})^{2\sum  l_i u_i}$.
\v8

 Let $(M, \dot{s})$ be a  spin  3-cobordism
from  $(\partial_-M,P_-, \sigma_-)$ to $(\partial_+M, P_+,
\sigma_+)$. The boundary of  $ M$ is  parametrized and
 $\dot{s}\in H^1(M,\partial M)$.
We assume that $\partial_-M$ (resp. $\partial_+M$) has $n_-$ (resp. $n_+$)
connected components.
We connect  the  marked points of  $\partial_- M$ (resp.
  $\partial_+ M$) by 
the   trivalent graph  $F_{n_-}$ (resp. a mirror image of $F_{n_+}$)
in $M$.
Let us  glue (along the parametrization)
to each connected component of
$ \partial_- M$ of genus $g$
a tubular neighborhood of the graph $\hat{G}^g$,
 containing the graph itself inside. 
The 1-vertex of $\hat{G}^g$ 
is  glued to the marked point. 
Analogously, to each connected component of
 $ \partial_+ M$ of genus $g$
we glue a tubular neighborhood of a mirror image of $\hat{G}^g$
(with respect to the plane orthogonal to that of the picture)
containing the graph itself  inside.
 The result is a closed 3-manifold $\tilde{M}$
with two closed 3-valent graphs  $\hat{G}^+$
and $\hat{G}^-$ inside.


We denote by $s$ the spin$^d$ structure on $M$ given by 
the homotopy class of $\tilde{\sigma}+\dot{s}$.
Let 
 $u_i$,
$2\leq i\leq n_+$, (resp. $u^\prime_i$, $2\leq i\leq n_-$)
be the number  associated by  $\dot{s}$
to the cycle in $M/\partial M$ obtained  
by identifying  the first and $i$th marked points of
$F_{n_+}$ (resp. $F_{n_-}$).

Let $\tilde{M}=S^3(L)$.
Let $\tilde{c}=(c_1,...,z_p)$
 be the solution of (\ref{re}) corresponding to $s$. 
Choose a special coloring $e$ (resp. $e^\prime$)
of $\hat{G}^+$ (resp. $\hat{G}^-$). Their
  grading on the circles is determined
 by $\{z_i\}$.
We define
$$
\tau_{e e^\prime}(M,\dot{s})
=\Delta^{-\s(L)} \eta^{-\chi(\partial_+ M)/2}
\sqrt{\langle e\rangle
\langle e^\prime\rangle}\;
u_e u^\prime_{e^\prime} \langle
L(\omega_{c_1},...,\omega_{c_m})\cup \hat{G}^+_e \cup \hat{G}^-_{e^\prime}\,
\rangle 
$$
and interpret it
 as an $(e, e^\prime)$-coordinate of the operator $\tau(M,\dot{s})$ 
from  the vector space
spanned by the special colorings of $\hat{G}^-$
to the vector space
spanned by the special colorings of $\hat{G}^+$.

The operator  $\tau(M,\dot{s})$
is an invariant of the spin 3-cobordism
$(M, \dot{s})$ with parametrized boundary. This is because,
it is an isotopy invariant of the graphs $\hat{G}^+$ and $\hat{G}^-$ and
it does not change under refined Kirby moves in $\tilde{M}$.  

We set  $\hat{G}=\hat{ G}^+\cup \hat{G}^-$. We call $L\cup \hat{G}$
the graph representing $M$, because
 $M$ can be reconstructed from it
(see [T, p.172]).

\begin{th}(Gluing property with anomaly)
If the 
spin$^d$ 3-cobordism
 $(M,\dot{s})$ is obtained from  $(M_1, \dot{s}_1)$
and $(M_2,\dot{ s}_2)$  by gluing along a diffeomorphism $f :\partial_+ M_1 
\to \partial_-  M_2 $  which preserves the relative spin$^d$ structures  and 
commutes with parametrizations, then 
\be \label{glu1}
 \tau_{{e}e^{\prime\prime} }(M,\dot{s})=
k\sum_{e^\prime}
\; \tau_{{e}e^{\prime} }(M_2,\dot{s}_2)\;
\tau_{{e^\prime} e^{\prime \prime} }(M_1, \dot{s}_1)\, ,\ee
where $k=\Delta^{\sigma(L)-\sigma(L_1)-\sigma(L_2)}$ is an anomaly factor
and $L$, $L_1$ and $L_2$ are the surgery links of $\tilde{M}$, 
$\tilde{M}_1$ and 
$\tilde{M}_2$
respectively.
 \end{th}    

\noindent
{\bf Remark.} To avoid the anomaly, 
we should supply
cobordisms  with so-called $p_1$-structures 
or 2-framings (see [BM] for more 
details). 
\v8\noindent
{\bf Proof:}
Let us suppose that $\partial_+ M_1$ has $n$ connected components.
We put the graph representing $M_2$  on  top of the graph 
representing  $M_1$
and introduce a
$0$-colored line connecting  $F_{n_\pm}$-lines of these graphs. Then 
we get
$$
\langle\mu\rangle \psdiag{7}{18}{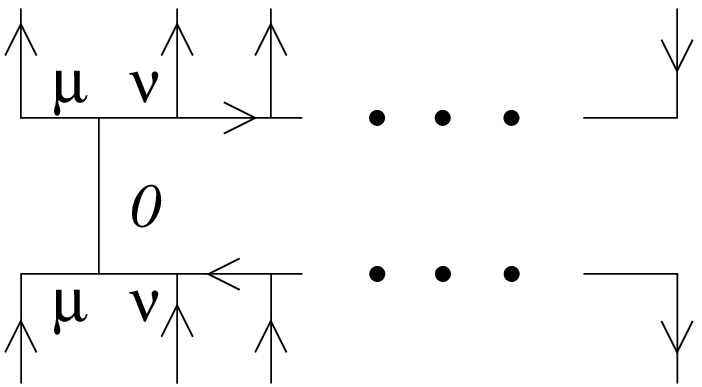}=\sum_{\l\mu}
 \langle \l \rangle\langle \mu \rangle
\psdiag{7}{18}{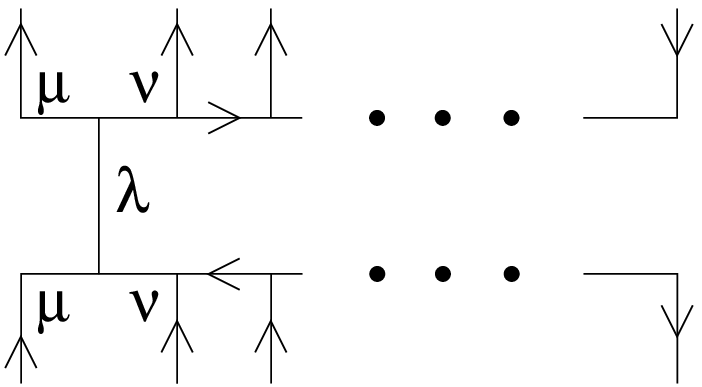}=\sum_{\mu}
\langle \mu \rangle\psdiag{7}{18}{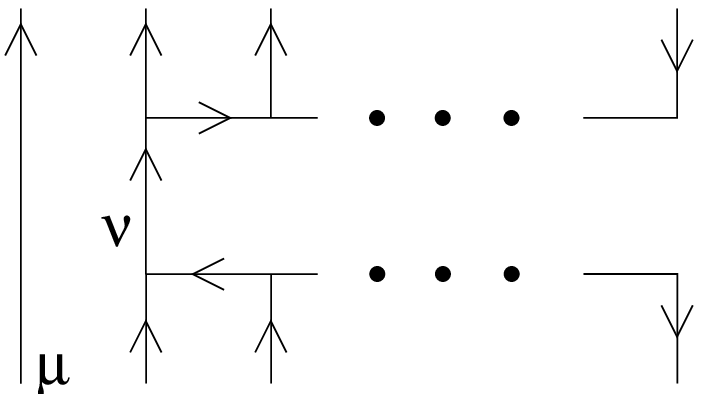}.$$ 
In the second equality we have used the fact that
for $\l\neq 0$ the Homfly polynomial of the colored graph is zero. 
The sum over all kinds of $(\l,\mu,\nu)$-vertices is assumed.
In (\ref{glu1}) the sum over all $\nu$ of the form
  $(1^N)^{\otimes k}\otimes
K^{l}$  is taken with $(a^N s)^{2 k (u_1+u^\prime_1)}(a^K s^{-1})^{2 l  (u_1
+u^\prime_1)}$ as 
coefficients, where $u_i$ (resp. $u^\prime_i$) is  assigned to the 
$i$-th $F_{n_+}$-line 
by $\dot{s}_1$ (resp. to the $i$th $F_{n_-}$-line by  $\dot{s_2}$)
and $u_1=u^\prime_1=0$ by  construction.
 This is equivalent to introducing a small
$\omega_0$-colored circle
around the  $\nu$-colored  line   and  allowing
 $\nu$ to run over $\Gamma_{N,K}$.
Continuing
 this procedure we will replace the figure drawn above by
  $n$ vertical strands, where
the
$i$th  strand ($2\leq i\leq n$)
is   linked with a small $\omega_{u_i+u^\prime_i}$-colored
  circle.
 After that,
 the sum over
all   colors of the remaining lines
 should be taken. Applying fusion rules 
again,  we get a graph representing $(M,s)$ (compare [T,p.177]). 
$\hfill\Box$
\v8

In $\cal S$ the identity morphism on $(\Sigma, P, \sigma)$
 is given by the cylinder $(\Sigma\times I, \dot{\sigma})$,
where $\dot{\sigma}$ is the trivial extension of $\sigma$. 
We define  $V (\Sigma, \sigma)$ to be the image of the projector 
 $\tau (\Sigma\times I, \dot{\sigma})$ associated to the cylinder.


The operator  $ \tau(M,\dot{s}):
V(\partial_-M, \sigma_-) \to V(\partial_+ M, \sigma_+)$
 defines
the  spin RT TQFT. 

\v8\noindent {\bf Remark.} 
In the spin TQFT
 the vector space associated with a non-connected surface 
with structure is not
equal to the tensor product of vector spaces associated with
connected components. 
 Therefore,
the  operators strongly depend on the cobordism structure of a
 given 3-manifold.
For example, the operators
$\tau(\Sigma\times I, \dot{s}):V(\Sigma,\sigma)\to V(\Sigma,\sigma)$
are equal for all extensions $\dot{s}$ of $\sigma$. But the operators
$\tau(\Sigma\times I, \dot{s}):V(\emptyset)\to V((-\Sigma,\sigma)\amalg
(\Sigma,
\sigma))$ distinguish $\dot{s}$.

\v8
\noindent
{\bf Weak spin RT TQFT.} 
Let us replace $\cal S$ with a weaker category,
where the objects are surfaces with spin$^d$  structure and
any  3-cobordism $M$ from $(\partial_- M,\sigma_-)$ to
$(\partial_+ M,\sigma_+)$
 is  supplied with $s\in {\rm Spin}^d(M)$, such that
$s|_{\partial_\pm M} =\sigma_\pm$.
Strictly speaking, it is not a category,
because
the spin$^d$ structure on the composition of such  cobordisms  
along $(\Sigma,\sigma)$ is  uniquely defined only 
 {\it  if $\Sigma$ is connected}. In order to get a category we should
allow cobordisms with a 'superposition' 
(or collection) of spin structures.

To define
the  invariant $\tau(M,s)$ we
  only need to  replace  $\hat{G}^g$ 
with $G^g$, depicted below, in the previous construction.  
\begin{center}
\mbox{\epsfysize=1.7cm\epsffile{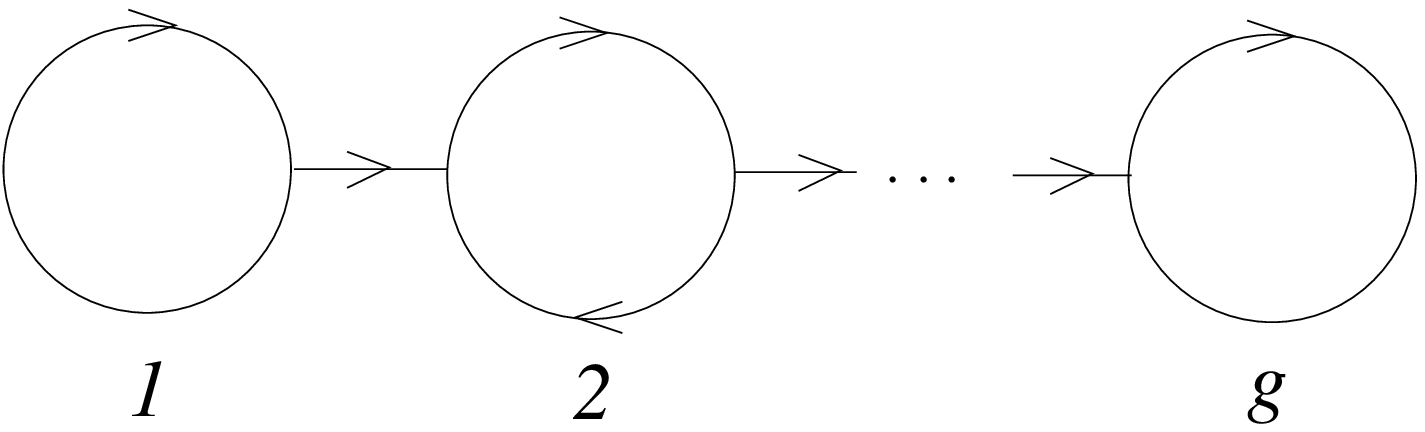}}
\\Figure 3 {\it The graph $G^g$.}
\end{center}
We denote by $G={ G}^+\cup {G}^-$ the resulting graph in $\tilde{M}$. 
The set of special colorings of $G$
 is a subset of the special colorings of $\hat{ G}$
consisting of colorings which are zero on $F_{n_+}\cup F_{n_-}$. 
The resulting
TQFT we shall call a weak spin RT TQFT.

In this TQFT, the vector space associated to 
the disjoint union of surfaces
is equal to the tensor product of spaces assigned to each of them. 
But we have 
 a weak form of the gluing property (compare Theorem 4 in [B]).

\begin{th}
If the 
 3-cobordism
 $(M,s)$ is obtained from  $(M_1, s_1)$
and $(M_2,{ s}_2)$  by gluing along a diffeomorphism $f :\partial_+ M_1 
\to \partial_-  M_2 $  which preserves spin$^d$ structures  and 
commutes with parametrizations, then 
\be 
\sum_s \tau_{{e}e^{\prime\prime} }(M,s)=
k\sum_{e^\prime}
\; \tau_{{e}e^{\prime} }(M_2,{s}_2)\;
\tau_{{e^\prime} e^{\prime \prime} }(M_1, {s}_1)\, ,\ee
where the sum is taken over all $s$ such that $s|_{M_i}=s_i$, $i=1,2$.
 \end{th}    

{\bf Remark.} The weak spin TQFT is the  `zero graded part'
of the spin TQFT  constructed in [BM]. The grading 
given by Theorem 11.2 in [BM] corresponds here to the 
orthogonal decomposition
of $V(\Sigma, \sigma)$ into subspaces generated by  colorings 
fixed on $F_{n_+}\cup F_{n_-}$.

\v8
For a 3-cobordism $M$,
we  define the vector $\tau(M)\in V(\partial M)$
 by its coordinates 
$$
\tau_{e e^\prime}(M)
=\Delta^{-\s(L)} \eta^{-\chi(\partial_+ M)/2}\sqrt{\langle e\rangle
\langle e^\prime\rangle}\;
 \langle G^+_e\cup
L(\omega,...,\omega)\cup G^-_{e^\prime}\rangle\, 
$$
in the basis of $V(\partial M)$ given by admissible colorings of $G$. 
The pair $(\tau(M), V(\partial M))$
defines  the unspun  RT TQFT.
Note that the number of admissible colorings of $G$
(given by Verlinde formula) 
 coincides with the dimension of  $V(\partial M)$
(see [L1]).
Analogous to Lemma 4.2 in [Bl], we 
can prove the  `transfer theorem', which identifies the unspun  theory
with the sum of weak  spin  TQFT's:
$$\tau(M)=\sum_{s\in {\rm Spin}^d (M)} \tau(M,s)$$

\v8\noindent
{\bf Refined TV  TQFT.}
Let us define a new cobordism
 category, where  an object is 
 ($\Sigma, P,\sigma,\tilde{ h}$)  
with  $\tilde{h}\in H^1(\Sigma, P)$ and
the structure on 3-cobordisms extends the one given  on the boundary.

 Let  $(M,\dot{s},\dot{h})$ be such a  3-cobordism
with  parametrized boundary.
Here $\dot{s}\in H^1(M,\partial M)$
determines the extension of
the relative spin$^d$ structure $\sigma$
 on $\partial M$ and
likewise,
 $\dot{h}\in H^1(M,\partial M)$ defines an extension of 
$\tilde{h}\in H^1(\partial M, P)$
to $M$.

We construct $(\tilde{M}, \hat{G})$
as in the spin RT TQFT. 
 Choose a handle decomposition of $\tilde{M}$
in such a way that $\hat{G}\subset H$.
Here $H$  is as before the union of 0- and 1-handles.

\v8 
The chain-mail graph  for  $(\tilde{M},\hat{G})$ 
is  the image under the embedding $j:H\hookrightarrow
 S^3$ of the graph consisting 
of 
\begin{itemize}
\item  a copy $\hat{G}^1$ of the graph $\hat{G}$ in the interior
 of $H$;
\item attaching curves of 2-handles pushed slightly into $H$;
\item a copy $\hat{G}^2$ of $\hat{G}$ on $\partial H$;
\item meridian curves of 1-handles pushed slightly into $S^3-H$.
\end{itemize}

The convention for the framing is the same as before. Denote by $A=\{a_1,...,
a_p\}$ the set of
circles of $\hat{G}$ and by $B$ the set of its meridians.
Let $u$ be the union of the sets $u$ and $u^\prime$
used in the spin RT TQFT.
Analogously, $t=\{0,t_2,...,t_{n_+},0,...,0\}\cup$
$\{0,t^\prime_2,...,t^\prime_{n_-},0,...,0\}$, where
 $t_i$ (resp. $t^\prime_i$)
is  the number
associated by  $\dot{h}$ to the cycle in $M/\partial M$
 obtained by identifying 
the first and $i$th marked points of
  $F_{n_+}$ 
(resp.  $F_{n_-}$). We denote by 
 $h\in H^1(M)$  the cohomology class determined by $\dot{h}$ and 
the cohomology class
on the boundary.

Then $s|_H-s_0|_{j(H)}$ assigns the numbers $\{x_1,...,x_a\}$ to the
 1-handles 
and the numbers $\{w_1,...,w_p\}$ to the elements of $B$.
Choose $y=\sum_i y_i\var_i$ representing $D(h)\in H_2(M,\partial M)$.
Then $\partial y=\sum_i v_i a_i$.
Choose a special  coloring $f$ (resp. $e$)
of $\hat{G}$, such that the grading of the colors on its $i$th
circle is equal to $-w_i$ (resp. $w_i-v_i$).

We set
\be\label{d}
Z_{ef}(M,\dot{s},\dot{h})=
(d\eta)^{d_0+d_3-2} \eta^{-\chi(\partial M)/2}
\sqrt{ \langle e\rangle\langle f\rangle}
u_f u_e t_e
\langle \hat{G}^2_f\cup R(\omega_{x_1}, ...,\omega_{y_b})\cup 
\hat{G}^1_e\rangle .\ee
We interpret $Z_{ef}(M,\dot{s},\dot{h})$ as an
 $(e,f)$-coordinate of the vector $Z(M,\dot{s},\dot{h})$ 
in  the space spanned  by special  colorings of the
graph $\hat{G}\cup \hat{G}$.

\begin{th}
$Z(M,\dot{s},\dot{h})$ is an invariant of the  3-cobordism
$(M,\dot{s},\dot{h})$ with
 parametrized  boundary. 
\end{th}

\noindent{\bf Proof:} By definition, $Z(M,\dot{s},\dot{h})$
 is an isotopy invariant of 
$\hat{G}^1\cup \hat{ G}^2$. The rest of the proof is analogous to the proof of Theorem 1.
Note that the number of
 lines in each 1-handle is $0$ mod $d$. Therefore the unknotting und 
reframing moves can be performed analogously.

Only the  births and deaths of 1-2-handle pairs require modifications.
It can happen that a birth of such a pair introduces
a $0$-framed $\omega_k$-colored  ($k=0,d/2$) unknot
 linked with 3-strands,
as depicted below:
\begin{center}
\mbox{\epsfysize=1.7cm\epsffile{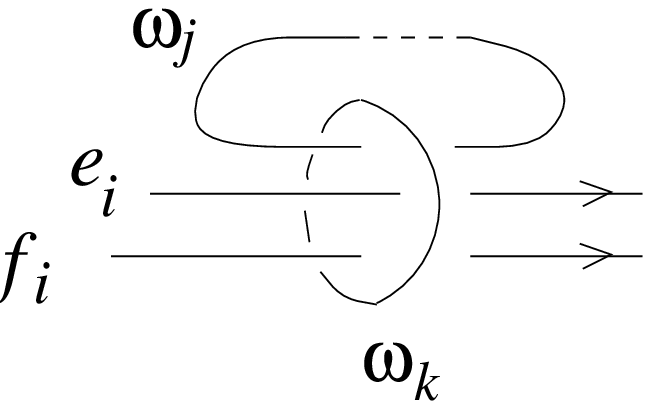}}
\end{center}

\noindent
We use the fusion rules to replace these strands by one.
Then applying 
 the graded
 killing property we get

$$\sum_{\l,\nu,\mu}\langle \l\rangle \langle \mu\rangle
\sum_{\alpha, \beta}
 \psdiag{8}{21}{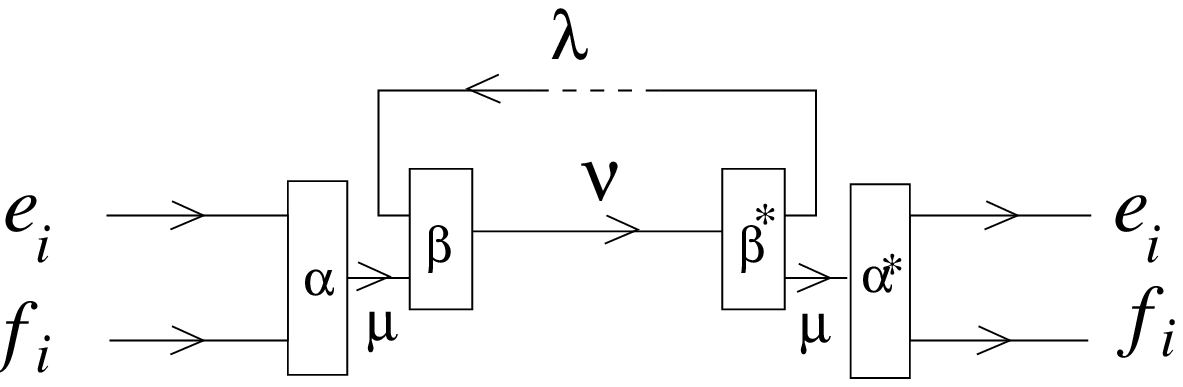},$$

\noindent 
 where we sum over all  $\nu$ of the form
 $(1^N)^{\otimes k}\otimes K^l$
and all $\l$ such that $|\l|=j$ mod $d$.  Note that $\langle \nu
\rangle=1$. Let us  apply (\ref{ff}) to the $\mu$-colored line.
After that, the sum $\sum_{\beta,\l}\langle \l\rangle=
\sum_{\l} N^{\nu}_{\l\mu}\langle \l\rangle=\langle \mu\rangle$
factorizes and
using  (\ref{f}) we can  delete the 1-2-handle pair.
$\hfill\Box$
\v8

\begin{th}\label{t}
For a  3-cobordism  $(M,\dot{s},\dot{h})$, 
$$Z(M,\dot{s},\dot{h})=\tau(-M,\dot{s})\otimes \tau(M,\dot{s}+\dot{h})$$
where  $\dot{s}+\dot{h}\in H^1(M,\partial M)$ 
 is the extension of  $\sigma +\tilde{h}$ to $M$.
\end{th}

The proof is analogous to the proof of Theorem 5. The difference is that
 the  handlebodies in  
(\ref{dec}) contain a copy of $\hat{G}$. 
\v8

Theorems \ref{t} and 7
 provide the gluing property (without anomaly)
for the invariant
 $Z(M,\dot{s},\dot{h})$. This completes 
the construction of the refined TV TQFT.

\v8
\noindent {\bf TV TQFT.}
Consider a 3-cobordism $M$ with parametrized boundary
 $\partial M=-\partial_-M\cup \partial_+ M$.
We construct $(\tilde{M}, G)$ 
as in the weak spin TQFT. 
 The admissible colorings of $G\cup G$ provide a  basis 
of  the vector space $V_{\partial M }$ associated with $\partial M$.

 Then  the vector
$Z(M)\in V_{\partial M}$ with coordinates
\be
Z_{e f}(M)=
\eta^{d_0+d_3-2} \eta^{-\chi(\partial_+ M)}
\sqrt {\langle e\rangle\langle f\rangle}
\langle G^2_f\cup R(\omega, ...,\omega)\cup G^1_e\rangle \ee
is an invariant of $M$ (by fogetting about the grading in the proof 
of Theorem 7). 
In fact,
$Z(M)$ is equal 
 to the invariant $\tilde{ Z}(M)$ defined in [BD1].
 This identifies the pair
 ($Z(M), V_{\partial M}$) with the Turaev-Viro TQFT.

We recall that $\tilde{Z}(M)$ is defined as the 
Turaev-Viro state sum operator of $M$ with 
 fixed triangulation of the boundary (given by  two copies of the dual 
graph to $G^g$ for each connected component of $\partial M$ of genus $g$).
The equality of $Z(M)$ and  $\tilde{Z}(M)$ can be shown
(in the spirit of Theorem 3.9 in [R1]) 
as follows:
Choose the dual triangulation of $\tilde{M}$ as  handle decomposition.
Using fusion rules and 
the killing property for 1-handle curves,
 we can split the graph  $G^2\cup R \cup G^1$
into parts sitting in  0-handles. 
This associates 6j-symbols to 0-handles with no 3-vertices of the graph
 inside and   products  of
6j-symbols to the others. 
 The definition of $\tilde {Z}(M)$
can then  be reconstructed  term-by-term.
 (The details will be omitted.)

The operator associated with a 3-cobordism $(M,s,h)$  by the weak refined
TV TQFT is denoted
 by $Z(M,s,h)$.

\begin{cor}
For a  3-cobordism  $(M,s,h)$, 
$$Z(M,s,h)=\tau(-M,s)\otimes \tau(M,s+h).$$
\end{cor}

\begin{cor}
The Turaev-Viro operator
invariant of a 3-cobordism $M$ splits into a sum of 
weak refined invariants, i.e.
$Z(M)=\sum_{s,h} Z(M,s,h)$.
\end{cor}



Finally,
we note that an  explicit  calculation
of a Homfly polynomial of a colored graph
requires
 the  knowledge of 
 6j-symbols  which   have apparently not yet been
determined for $N>2$. 

\v8\v8
\noindent
{\sl Mathematisches Institut, Universit\"at Bern, Sidlerstr. 5, 
CH-3012 Bern \\ 
e-mail: {\tt beliak$\char'100$math-stat.unibe.ch}}


\end{document}